\documentclass[12pt]{amsart}
\usepackage{mathptmx}
\usepackage{amsmath}
\usepackage{amssymb}
\usepackage{array}
\usepackage{geometry}
\usepackage[bookmarks=true,colorlinks=true, pdfstartview=FitV, linkcolor=black, citecolor=blue, urlcolor=black]{hyperref}
\usepackage{color}
\definecolor{DarkRed}{rgb}{0.55,.00,0.2}
\definecolor{DarkGrey}{rgb}{0.35,.35,0.35}

\theoremstyle{definition}

\theoremstyle{remark}

\numberwithin{equation}{section}

\begin{document}

\title{On the  half-Hartley transform,  its iteration and compositions with  Fourier transforms}

\author{S. Yakubovich}
\address{Department of Mathematics, Fac. Sciences of University of Porto,Rua do Campo Alegre,  687; 4169-007 Porto (Portugal)}
\email{ syakubov@fc.up.pt}



\maketitle

\markboth{\rm \centerline{ S.  YAKUBOVICH}}{}
\markright{\rm \centerline{ COMPOSITIONS OF THE HARTLEY AND FOURIER TRANSFORMS }}

\begin{abstract} Employing the generalized Parseval equality for the Mellin transform and elementary trigonometric formulas,
the iterated Hartley transform on the nonnegative  half-axis (the iterated half-Hartley transform) is investigated in $L_2$.
 Mapping and inversion properties are discussed, its relationship with the iterated Stieltjes transform is established.
 Various compositions with the Fourier cosine and sine transforms are obtained. The results are applied to  the uniqueness and universality
 of the closed form solutions  for certain  new singular integral and integro-functional equations.
\bigskip
\end{abstract}

\section{Introduction and auxiliary results}

The familiar reciprocal pair of the Hartley transforms \cite{Brace}
$$(\mathcal{H} f)(x) = {1\over \sqrt {2\pi} } \int_{-\infty}^\infty [\cos(xt)+ \sin(xt) ] f(t) dt, \ x \in \mathbb{R},\eqno(1.1)$$
$$f(x)= {1\over \sqrt {2\pi} } \int_{-\infty}^\infty [\cos(xt)+ \sin(xt) ] (\mathcal{H} f)(t) dt\eqno(1.2)$$
is well-known  in connection with various  applications in mathematical physics.  Mapping  and inversion properties of these transforms in $L_2$ as well as their multidimensional analogs were investigated, for instance, in \cite{vu}, \cite{haiw}, \cite{luch}. These operators were treated as the so-called bilateral Watson transforms,  and in some sense they are related to the Fourier cosine and Fourier sine transforms
$$(F_c f)(x) = \sqrt {{2\over \pi} } \int_{0}^\infty \cos(xt)  f(t) dt, \ x \in \mathbb{R}_+,\eqno(1.3)$$
$$(F_s f)(x) = \sqrt { {2\over \pi} } \int_{0}^\infty  \sin(xt)  f(t) dt,\   x \in \mathbb{R}_+.\eqno(1.4)$$
Recently \cite{hart}, the author investigated  the Hartley transform (1.1) with the integration over $\mathbb{R}_+$ (the half-Hartley transform)
$$(\mathcal{H}_+ f)(x) = \sqrt { {2\over \pi} } \int_{0}^\infty [\cos(xt)+ \sin(xt) ] f(t) dt, \ x \in \mathbb{R}_+\eqno(1.5)$$
and proved and analog of the Plancherel theorem, establishing  its reciprocal inverse operator  in $L_2(\mathbb{R}_+)$
$$f(x)=   \sqrt { {2\over \pi} }  \int_0^\infty \left[ \sin( xt)  \ S(xt)+ \cos (xt) \  C(xt)\right] (\mathcal{H}_+ f)(t) dt,\eqno(1.6)$$
where $S(x),\ C(x)$ are Fresnel sin- and cosine- integrals, respectively,
$$S(x)=  \sqrt { {2\over \pi} } \int_0^{\sqrt x}  \sin(t^2) dt, \quad   C(x)=  \sqrt { {2\over \pi} } \int_0^{\sqrt x}  \cos (t^2) dt.$$

Our goal here is to examine the iteration of the operator (1.5)   $(\mathcal{H}^2_+ f)(x) \equiv (\mathcal{H}_+ \mathcal{H}_+ f)(x)$, which will be called the iterated half-Hartley transform and investigate its compositions in $L_2$ with the Fourier transforms (1.3),  (1.4) of the form:  $\mathcal{H}_+ F_c ,\  \mathcal{H}_+ F_s ,\  \mathcal{H}_+ F_c F_s, \  \mathcal{H}^2_+ F_c,\   \mathcal{H}^2_+ F_s,\   \mathcal{H}^2_+ F_c F_s$.  The corresponding integral representations of these compositions will be established in $L_2$ and their boundedness and invertibility will be proved. Moreover, we will apply these results to establish the uniqueness of solutions in the closed form for  the corresponding second kind  singular integral and integro-functional equations.

Our  natural approach is based on  the $L_2$-theory of the Mellin transform \cite{tit}
$$(\mathcal {M} f)(s)= f^*(s)= \int_0^\infty f(t) t^{s-1}dt, \ s \in
\sigma  =\{ s \in \mathbb{C}, \   s={1\over 2}  +i\tau\},\eqno(1.7)$$
where the integral  is convergent in the mean square sense with respect to the norm in $L_2(\sigma)$. Reciprocally,  the inversion formula takes place
$$f(x)= {1\over 2\pi i}\int_{\sigma} f^*(s) x^{-s} ds,\ x >0\eqno(1.8)$$
with the convergence of the integral in the mean square sense with respect to the norm in $L_2(\mathbb{R}_+)$.  Furthermore, for any $f_1, f_2  \in L_2(\mathbb{R}_+)$ the generalized Parseval identity holds
$$\int_0^\infty f_1\left(xt\right) f_2(t) dt  = {1\over 2\pi i}\int_{\sigma}  f_1^*(s)f_2^*(1-s) x^{-s}
ds, \ x >0\eqno(1.9)$$
with  Parseval's equality of squares of $L_2$- norms
$$\int_0^\infty |f(x)|^2 dx = {1\over 2\pi}   \int_{-\infty}^{\infty} \left|f^*\left({1\over 2}  + i\tau\right)\right|^2 d\tau.\eqno(1.10)$$
Finally in this section, we exhibit  the known formulas \cite{tit}
$$\int_0^\infty {\sin t \over t} t^{s-1} dt = \frac{\Gamma(s)}{1-s} \cos\left({\pi s\over 2}\right), \ s \in \sigma,\eqno(1.11)$$
$$\int_0^\infty {1-\cos t \over t} t^{s-1} dt = \frac{\Gamma(s)}{1-s} \sin \left({\pi s\over 2}\right), \ s \in \sigma,\eqno(1.12)$$
where $\Gamma(s)$ is the Euler gamma-function,  which will be used in the sequel.

\section{ Plancherel's type theorems}

It is widely known \cite{tit} via the classical Plancherel theorem in $L_2(\mathbb{R}_+)$ that the Fourier cosine and Fourier sine transformations extend to bounded invertible and isometric mappings, having the properties $F_c^2 =I, \  F_s^2= I$, where $I$ is the identity operator.   We begin, demonstrating our method on the simple composition $F_cF_s$ of operators (1.3), (1.4).  Precisely, it drives us to the Plancherel theorem for the Hilbert transform \cite{tit}.

{\bf Theorem 1}. {\it  The composition $F(x)= (F_c F_s f)(x)$  extends to a bounded invertible and isometric map
$F:  L_2(\mathbb{R}_+) \to L_2(\mathbb{R}_+)$ and can be written in  the form of the Hilbert transform in $L_2$
$$F(x) =  {2\over \pi}  PV  \int_{0}^\infty {t f(t) \over t^2- x^2}  dt, \quad  x \in \mathbb{R}_+. \eqno(2.1)$$
Reciprocally,
$$f(x) =  {2\over \pi}  PV \int_{0}^\infty {x \  F(t) \over x^2- t^2}  dt, \quad  x \in \mathbb{R}_+ \eqno(2.2)$$
and this map is isometric, i.e. $||F|| = ||f|| $ for all $f \in  L_2(\mathbb{R}_+)$. }

\begin{proof}   Let $f$ belong to the space $C^{(2)}_c(\mathbb{R}_+)$  of  continuously differentiable functions
of compact support, which is dense in $ L_2(\mathbb{R}_+)$.  Then integrating by parts  in (1.3), (1.4),  (1.7),   we find that $(F_c f)(x)= O(x^{-2} ),\  (F_s  f)(x)= O(x^{-2} ),\ x \to \infty$ and  $s^2 f^*(s)$ is bounded on $\sigma$.  Therefore $(F_c f)(x), (F_s f)(x) \in   L_2(\mathbb{R}_+) \cap  L_1(\mathbb{R}_+), \   f^*(s) \in L_2(\sigma) \cap L_1(\sigma) $.   Hence minding equalities (1.11), (1.12),  the generalized Parseval equality (1.9) and the supplement formula for the gamma-function, we derive for all   $x >0$ (cf. \cite{tit}, Section 8.4)
$$F(x)= {1\over 2\pi i}\int_{\sigma}  F^*(s) x^{-s} ds  = {2\over \pi} {1\over 2\pi i} \int_{\sigma}  \Gamma(s) \Gamma(1-s) \cos^2 \left({\pi s\over 2}\right) f^*(s) x^{-s} ds$$$$
=   {1\over 2\pi i}\int_{\sigma}  f^*(s) \cot \left({\pi s\over 2}\right) x^{-s} ds =   {2\over \pi}  PV \int_{0}^\infty {t f(t) \over t^2- x^2}  dt .$$
Hence, reciprocally via (1.8),  we obtain
$$f(x)= {1\over 2\pi i}\int_{\sigma}  f^*(s) x^{-s} ds =  {1\over 2\pi i}\int_{\sigma}  F^*(s)
\tan \left({\pi s\over 2}\right) x^{-s} ds $$$$=  {1\over 2\pi i}\int_{\sigma}  F^*(1-s)
\cot \left({\pi s\over 2}\right) x^{s-1} ds = {2\over \pi x}  PV \int_{0}^\infty { F(t) \over 1- (t/x)^2}  dt =  {2\over \pi }  PV \int_{0}^\infty {x F(t) \over x^2- t^2}  dt .$$
Thus  we proved (2.1),  (2.2)  for  any $f \in C^{(2)}_c(\mathbb{R}_+)$.   Further, since $C^{(2)}_c(\mathbb{R}_+)$ is dense in $ L_2(\mathbb{R}_+)$, there is a unique extension of $F$ as an invertible continuous map $F:  L_2(\mathbb{R}_+) \to L_2(\mathbb{R}_+)$. Clearly, it is isometric by virtue of the Plancherel theorem for Fourier cosine and Fourier sine transforms.

\end{proof}

Extending  this approach, we  prove  the Plancherel type theorem for the iterated half-Hartley transform $\mathcal{H}^2_+$.   Indeed, we have

{\bf Theorem 2}. {\it  The iterated half-Hartley transform  extends to a bounded invertible map
$\mathcal{H}^2_+:  L_2(\mathbb{R}_+) \to L_2(\mathbb{R}_+)$ and the following reciprocal formulas hold
$$(\mathcal{H}^2_+ f)(x) = 2 f(x) +  {2\over \pi} \int_0^\infty \frac {f(t)} {x+t} dt,\quad  x > 0,\eqno(2.3)$$
$$f(x)=   {1\over 2} (\mathcal{H}^2_+ f) (x) -    {1 \over \pi^2}  \int_{0}^\infty  \frac{ \sqrt{xt} \  [ \log x-  \log t] } { x^2-  t^2 }   (\mathcal{H}^2_+ f) (t)  dt.\eqno(2.4)$$
Moreover,  the  norm inequalities take place}
$$ \left|\left| f \right|\right|  \le \left|\left|  \mathcal{H}^2_+ f \right|\right| \le
  8 \   \left|\left| f \right|\right| .\eqno(2.5)$$

\begin{proof}  Assuming again $f \in C^{(2)}_c(\mathbb{R}_+)$ and taking into account (1.8), (1.9) and relation (8.4.2.5) in \cite{prud},   Vol. 3,  we derive in the same manner  the equalities
$$(\mathcal{H}^2_+ f)(x) = {1\over 2 \pi i}  {2\over \pi} \int_{\sigma}  \Gamma(s)\Gamma(1-s)\left(\cos \left({\pi s\over 2}\right) + \sin \left({\pi s\over 2}\right) \right)^2 f^*(s) \  x^{-s} ds$$
$$=   {1\over \pi i}  \int_{\sigma}  \frac {1+  \sin \left(\pi s \right) } { \sin \left(\pi s \right) }  f^*(s) \  x^{-s} ds = 2 f(x)+ {1\over 2 \pi i}  {2\over \pi}   \int_{\sigma}  \Gamma(s)\Gamma(1-s)  f^*(s)  \  x^{-s} ds$$
$$=  2 f(x) +  {2\over \pi} \int_0^\infty \frac {f(t)} {x+t} dt,\eqno(2.6)$$
which prove representation (2.3),  involving the classical Stieltjes transform \cite{luch},  \cite{tit}.  Conversely,  appealing to relation (8.4.6.11) in \cite{prud},  Vol. 3 and elementary properties of the Mellin transform, it gives

$$ f(x)=  {1\over 4 \pi i}  \int_{\sigma}  \frac { \sin \left(\pi s \right) }  { 1+  \sin \left(\pi s \right) }  (\mathcal{H}^2_+ f)^*(s) \  x^{-s} ds =  {1\over 2} (\mathcal{H}^2_+ f) (x)$$

$$ -   {1\over 8 \pi i}  \int_{\sigma}  \frac { (\mathcal{H}^2_+ f)^*(s)  }  {    \sin^2 \left(\pi (s +1/2)/2  \right)} \  x^{-s} ds =   {1\over 2} (\mathcal{H}^2_+ f)^*(x)$$

$$-   {1\over 2 \pi i} {1\over 4\pi^2}  \int_{\sigma}   \left[  \Gamma\left ({s\over 2} + {1\over 4} \right)\Gamma\left({3\over 4} - {s\over 2}\right)\right]^2  (\mathcal{H}^2_+ f)^*(s)   \  x^{-s} ds  $$

$$ =  {1\over 2} (\mathcal{H}^2_+ f) (x) -    {1 \over \pi^2}  \int_{0}^\infty  \frac{ \sqrt{xt} \  [ \log x-  \log t ]} { x^2-  t^2 }   (\mathcal{H}^2_+ f) (t)  dt, \eqno(2.7)$$
which proves (2.4), involving the iterated Stietltjes transform recently treated in \cite{martin}.   In order to establish inequalities (2.5), we call the Parseval equality (1.10) for the Mellin transform, which yields  (see (2.6), (2.7))
$$ \left|\left|  \mathcal{H}^2_+ f \right|\right|  =  2  \left( {1\over 2\pi} \int_{-\infty}^{\infty}  \frac {\left(1+  \cosh \left(\pi \tau \right) \right)^2} { \cosh^2 \left(\pi \tau \right) }  \left| f^*\left({1\over 2} +i\tau\right)\right|^2 d\tau\right)^{1/2}$$
$$= 4  \left( {1\over 2\pi} \int_{-\infty}^{\infty}  \frac { \cosh^4 \left(\pi \tau /2\right) } { \cosh^2 \left(\pi \tau \right) }  \left| f^*\left({1\over 2} +i\tau\right)\right|^2 d\tau\right)^{1/2}\le 8 \   ||f||,$$
and  on the other hand,
$$ ||f||=    {1\over 4 \sqrt {2\pi} }  \left( \int_{-\infty}^{\infty}  \frac  { \cosh^2 \left(\pi \tau \right) }    { \cosh^4 \left(\pi \tau /2\right) }\left|  (\mathcal{H}^2_+ f)^* \left({1\over 2} +i\tau\right)\right|^2 d\tau\right)^{1/2} \le  \left|\left|  \mathcal{H}^2_+ f \right|\right|.$$
Now the same argument of the denseness of  $C^{(2)}_c(\mathbb{R}_+)$ in $ L_2(\mathbb{R}_+)$ drives us to  a unique extension of $\mathcal{H}^2_+$ as an invertible continuous map $\mathcal{H}^2_+:  L_2(\mathbb{R}_+) \to L_2(\mathbb{R}_+)$.
\end{proof}

{\bf Remark 1}.  As we observe via the Schwarz inequality,  the convergence of integrals (2.3), (2.4) is pointwise.

Concerning the composition $\mathcal{H}_+ F_c$  it has

{\bf Theorem 3}. {\it  The composition $F(x)= \left(\mathcal{H}_+ F_c  f\right)(x) $  extends to a bounded invertible map
$F:  L_2(\mathbb{R}_+) \to L_2(\mathbb{R}_+)$ and  the following reciprocal formulas hold
$$F (x) =  f(x) +  {2\over \pi }  PV \int_{0}^\infty {x f(t) \over x^2- t^2}  dt,\quad  x > 0, \eqno(2.8)$$
$$f(x)=   {1\over \pi}  PV \int_0^\infty  {\sqrt { x t} \over t^2-  x^2 }  F(t)  dt .\eqno(2.9)$$
Moreover,  the  norm inequalities are valid}
$$ \left|\left| f \right|\right|  \le \left|\left|  F \right|\right| \le   2 \sqrt 2    \left|\left| f \right|\right| .\eqno(2.10)$$

\begin{proof}  For  $f \in C^{(2)}_c(\mathbb{R}_+)$ we obtain
$$F(x) = {1\over 2 \pi i}  {2\over \pi} \int_{\sigma}  \Gamma(s)\Gamma(1-s)\left(\cos \left({\pi s\over 2}\right) + \sin \left({\pi s\over 2}\right) \right) \sin \left({\pi s\over 2}\right) f^*(s) \  x^{-s} ds$$
$$=   {1\over 2 \pi i}  \int_{\sigma}  \left(1+  \tan \left({\pi s\over 2}\right) \right)  f^*(s) \  x^{-s} ds =  f(x)+ {2\over \pi }  PV  \int_{0}^\infty {x f(t) \over x^2- t^2}  dt,$$
which prove representation (2.8),  relating again to the classical Hilbert  transform in $L_2(\mathbb{R}_+)$.  The inverse operator (2.9) can be deduced via the equality

$$ f(x)=  {1\over 2 \pi i}  \int_{\sigma}  \frac { F^*(s) }  {  1+  \tan \left(\pi s/2 \right) }  \  x^{-s} ds$$

 $$=  {1\over 2\sqrt 2 \  \pi i}  {d\over dx} \int_{\sigma}  F^*(s) \frac{\Gamma\left ((s+ 1/2)/ 2 \right)\Gamma\left((3/2 - s)/2\right)}{\Gamma\left ((1+s)/2 \right)\Gamma\left((1 - s)/2\right)}    \  {x^{1-s}\over 1-s} ds,\eqno(2.11)$$
where the differentiation is allowed under the integral sign via the absolute and uniform convergence.  In the mean time,    owing to the residue theorem

$${1\over 2\pi i}  \int_{\sigma}  \frac{\Gamma\left ((s+ 1/2)/ 2 \right)\Gamma\left((3/2 - s)/2\right)}{\Gamma\left ((1+s)/2 \right)\Gamma\left((1 - s)/2\right)}    \  {x^{1-s}\over 1-s}  ds =  {\sqrt 2 \over \pi} \sum_{n=0}^\infty  {x^{2n+3/2}\over 2n + 3/2} =  {\sqrt 2  \over \pi} \int_0^x  {\sqrt y \  dy \over 1-y^2},\   0 < x < 1, $$
and

$$ {1\over 2\pi i}  \int_{\sigma}  \frac{\Gamma\left ((s+ 1/2)/ 2 \right)\Gamma\left((3/2 - s)/2\right)}{\Gamma\left ((1+s)/2 \right)\Gamma\left((1 - s)/2\right)}    \  {x^{1-s}\over 1-s}  ds =     {\sqrt {2 } \over \pi} \sum_{n=0}^\infty { x^{- 2n- 1/2}\over 2n+1/2} =  {\sqrt  2 \over \pi}  \int_x^\infty   {\sqrt y \  dy \over y^2-1},\    x > 1.$$
Therefore, returning to (2.11), differentiating with respect to $x$ under the integral sign and using (1.9), we write it in the form
$$f(x)=  {1\over \pi} PV \int_0^\infty  {\sqrt { x/ t} \over 1-  (x/t) ^2 }  F(t)  {dt\over t} $$
and come out with (2.9).  Finally, in a similar manner, we derive  inequalities (2.10). In fact,  we have
$$ ||f|| \le  \left|\left|  F \right|\right| =   2   \left( {1\over 2\pi} \int_{-\infty}^{\infty}  \frac { \cosh^2 \left(\pi \tau/2 \right) } { \cosh \left(\pi \tau \right) }  \left| f^*\left({1\over 2} +i\tau\right)\right|^2 d\tau\right)^{1/2} \le 2 \sqrt 2  || f||.$$
Hence extending  $F$ on the whole  $ L_2(\mathbb{R}_+)$ as an invertible continuous mapping, we complete the proof.
\end{proof}

The Plancherel theorem for the composition  $\mathcal{H}_+ F_s$ can be proved analogously,  and we leave   it without proof.

{\bf Theorem 4}. {\it  The composition $F(x)= \left(\mathcal{H}_+ F_s  f\right)(x) $  extends to a bounded invertible map
$F:  L_2(\mathbb{R}_+) \to L_2(\mathbb{R}_+)$ and  the following reciprocal formulas hold
$$F (x) =  f(x) +  {2\over \pi }  PV \int_{0}^\infty {t f(t) \over t^2- x^2}  dt,\quad  x > 0, \eqno(2.12)$$
$$f(x)=   {1\over \pi}  PV \int_0^\infty  {\sqrt { x t} \over x^2-  t^2 }  F(t)  dt .\eqno(2.13)$$
Moreover,  the  norm inequalities $(2.10)$ are valid}.

The case  $\mathcal{H}_+ F_c F_s$ can be treated with the use of Theorem 1.   Precisely, we state

{\bf Theorem 5}. {\it  The composition $F(x)= \left(\mathcal{H}_+ F_c F_s  f\right)(x) $  extends to a bounded invertible map  $F:  L_2(\mathbb{R}_+) \to L_2(\mathbb{R}_+)$, having   the integral representation
$$F (x) = \sqrt { {2\over \pi} } \int_{0}^\infty \left[ \sin(xt) - \cos(xt)+ {2\over \pi}  \sqrt {xt} \   S_{-1/2,1/2} (xt) \right] f(t) dt,\quad  x > 0, \eqno(2.14)$$
where $ S_{-1/2,1/2} (x)$ is the Lommel function and the integral converges in the mean square sense.    The inverse operator is written in the form of the integral

$$f(x)=  2  \sqrt{{2\over \pi}} \int_0^\infty  \  \left[ (1-  S(xt)) \sin( xt)  -   C(xt) \cos( xt)  \right]  F(t)dt,\eqno(2.15)$$
which converges in the mean square sense as well.  Moreover,  the  norm inequalities $(2.10)$ take place}.

\begin{proof}  In fact,  for  $f \in C^{(2)}_c(\mathbb{R}_+)$ we write via (1.9)
$$F(x) = {1\over 2 \pi i}  \left({2\over \pi}\right)^{3/2}  \int_{\sigma}  \Gamma^2(s)\Gamma(1-s)\left(\cos \left({\pi s\over 2}\right) + \sin \left({\pi s\over 2}\right) \right) \sin^2 \left({\pi s\over 2}\right) f^*(1- s) \  x^{-s} ds$$
$$=   {1\over 2 \pi i}  \sqrt {{2\over \pi}} \int_{\sigma}   \Gamma (s) \sin \left({\pi s\over 2}\right) \left(1+  \tan \left({\pi s\over 2}\right) \right)  f^*(1- s) \  x^{-s} ds=  (F_sf)(x)- (F_c f)(x)$$
$$+   {1\over 2 \pi i}  \sqrt {{2\over \pi}} \int_{\sigma}  \frac { \Gamma (s) } { \cos \left(\pi s/2\right) }  f^*(1- s) \  x^{-s} ds.\eqno(2.16)$$
Meanwhile, the latter integral can be calculated, appealing to   relations (8.4.2.5) in \cite{prud},  Vol. 3,   (8.4.23.1) in \cite{prud},  Vol. 3 and (2.16.3.14) in \cite{prud},  Vol. 2.   But first,  employing the supplement and duplication formulas  for the gamma-function,   we find the following inverse Mellin transform   written in terms of the Mellin type convolution with the modified Bessel function, namely
$$ {1\over 2 \pi i}   \int_{\sigma}  \frac { \Gamma (s) } { \cos \left(\pi s/2\right) }  \  x^{-s} ds =
{1\over  4 \pi^{5/2} i}  \int_{\sigma}   \Gamma^2 \left({1+s\over 2}\right)  \Gamma \left({s\over 2}\right)\Gamma\left({1-s\over 2}\right)  \ ( x/2)^{-s} ds $$
$$= {2\over   \pi}  \int_{0}^\infty   \frac{\sqrt t\   K_0(2\sqrt t )}{ \sqrt {1+ (x^2/ 4t)}} {dt\over t} = {2\over  \pi}  \int_{0}^\infty   \frac{ y \  K_0(y )}{ \sqrt {y^2 + x^2}} dy =  { 2  \sqrt x\over  \pi}  \   S_{-1/2,1/2} (x),\   x >0,$$
where $S_{\mu,\nu} (z) $ is the Lommel function  \cite{prud},  Vol. 3.  Therefore,  returning to (2.16) and calling  again the generalized Parseval equality (1.9),  we get

$${1\over 2 \pi i}  \sqrt {{2\over \pi}} \int_{\sigma}  \frac { \Gamma (s) } { \cos \left(\pi s/2\right) }  f^*(1- s) \  x^{-s} ds= {2 \sqrt 2\over \pi\sqrt \pi } \int_0^\infty \sqrt {xt} \   S_{-1/2,1/2} (xt) f(t) dt, \quad  x >0,$$
 and the latter integral is absolutely convergent for any $f \in L_2(\mathbb{R}_+)$.  Combining with (2.16), we come out with  representation (2.14) for the dense set  $C^{(2)}_c(\mathbb{R}_+)$ of  $L_2(\mathbb{R}_+).$  Moreover, the norm inequalities (2.10) follow immediately from (2.16).  In fact,  we have

$$ ||f|| \le ||F||=     {\sqrt 2\over  \pi } \left(  \int_{-\infty}^\infty   \left| \Gamma \left({1\over 2} +i\tau \right) \right|^2 \cosh^2(\pi\tau/2)  \left|f^*\left({1\over 2} + i\tau\right)\right|^2  d\tau \right)^{1/2} $$

$$=  \sqrt {{ 2\over  \pi }}  \left( \int_{-\infty}^\infty   \frac{ \cosh^2(\pi\tau/2) } {\cosh(\pi\tau)}  \left|f^*\left({1\over 2} + i\tau\right)\right|^2  d\tau \right)^{1/2} \le 2\sqrt 2 ||f||. $$
Concerning the inverse operator (2.15),  we call again (2.16) and reciprocal formulas (1.7), (1.8) of the Mellin transform.  It yields

$$ f(x)=   {1\over 2 \pi i}  \sqrt {{\pi\over 2}} \int_{\sigma}   \frac{ F^*(1- s)} { \Gamma (1-s) \cos \left(\pi s/ 2\right) \left(1+  \cot \left(\pi s/ 2\right) \right)}  \  x^{-s} ds$$

$$=  {1\over 4 \pi i}  \int_{\sigma}   \frac{\Gamma\left ((s+ 1)/ 2 \right)\Gamma\left ((s+ 1/2)/ 2 \right)\Gamma\left((3/2 - s)/2\right)}{\Gamma\left (s/2 \right)\Gamma^2\left(1 - s/2\right)}  F^*(1- s)  \  (x/2)^{-s} ds.\eqno(2.17)$$
Meanwhile, the  residue theorem and relation (7.14.4.6) in \cite{prud},  Vol. 3 lead us to the value of the integral with the ratio of gamma-functions in terms of Fresnel's integrals (see above)
$$ {1\over 2 \pi i}  \int_{\sigma}   \frac{\Gamma\left ((s+ 1)/ 2 \right)\Gamma\left ((s+ 1/2)/ 2 \right)\Gamma\left((3/2 - s)/2\right)}{\Gamma\left (s/2 \right)\Gamma^2\left(1 - s/2\right)}  (x/2) ^{-s} ds = 2  \sqrt{{2\over \pi}} \  \left[ (1-  S(x)) \sin x \right.$$$$\left.  -   C(x) \cos x \right] .$$
Hence, returning to (2.17), we easily come out with (2.15),  and after  extension of  $F$ on the whole  $ L_2(\mathbb{R}_+)$ as an invertible continuous mapping  complete the proof of Theorem 5.
\end{proof}

The Plancherel theorem for compositions  $ \mathcal{H}^2_+ F_c,\    \mathcal{H}^2_+ F_s$ is related to Theorem 2 and can be stated as follows

{\bf Theorem 6}. {\it  Compositions  $F(x)= \left(\mathcal{H}^2_+ F_c   f\right)(x),\  G(x)=  \mathcal{H}^2_+ F_s $  extend to  bounded invertible mappings   $F, G:  L_2(\mathbb{R}_+) \to L_2(\mathbb{R}_+)$, having    integral representations
$$F (x) = 2 \sqrt { {2\over \pi} } \int_{0}^\infty \left[  \cos(xt)+ {1\over \pi}
 \sqrt {xt} \   S_{-3/2, -1/2} (xt) \right] f(t) dt,\quad  x > 0, \eqno(2.18)$$

$$G (x) = 2 \sqrt { {2\over \pi} } \int_{0}^\infty \left[  \sin(xt)+ {1\over \pi}  \sqrt {xt} \
 S_{-1/2, 1/2} (xt) \right] f(t) dt,\quad  x > 0, \eqno(2.19)$$
where both  integrals converge in the mean square sense.    Inverse operators are  written, respectively,  in the form

$$f(x)=     \int_{0}^\infty  k_c(xt) F (t)  dt,\quad x > 0\eqno(2.20)$$
where
$$k_c(x)=  {\sqrt x \over  \pi}  \sum_{k=0}^\infty   {(-1)^k \ x^{2k} \over (3/2)_ {2k}} \left[  {2\over \pi} \psi( -1/2- 2k)  -  {2\over \pi} \log x   + 1 \right]$$
and  $\psi(x)$ is the psi-function,
$$f(x)=   \int_{0}^\infty k_s(xt)  G (t) dt,\eqno(2.21)$$
where
$$k_s(x)= {\sqrt x \over   \pi}   \sum_{k=0}^\infty  { (-1)^k \  x^{ 2k} \over (3/2)_ {2k}} \left[ 1-  {2\over \pi} \psi( -1/2- 2k)  +  {2\over \pi} \log x  \right]$$
and  integrals  converge in the mean square sense.   Moreover,  the  norm inequalities $(2.5)$ take place}
$$   \left|\left| f \right|\right|  \le \left\{ \begin{array}{c}  \left|\left|  F \right|\right| \\
\left|\left|  G \right|\right| \end{array} \right\} \le  8  \left|\left| f \right|\right|.\eqno(2.22)$$

\begin{proof}  Indeed,  for  $f \in C^{(2)}_c(\mathbb{R}_+)$ we use (2.6) to obtain
$$F(x)  = {1\over 2 \pi i}  \left({2\over \pi}\right)^{3/2} \int_{\sigma}  \Gamma(s)\Gamma^2(1-s)\left(\cos \left({\pi s\over 2}\right)
 + \sin \left({\pi s\over 2}\right) \right)^2 \sin \left({\pi s\over 2}\right) f^*(s) \  x^{s-1} ds$$
$$=   {1\over 2\pi i} \sqrt {{2\over \pi}} \int_{\sigma}  \frac {1+  \sin \left(\pi s \right) } { \sin \left(\pi s/2 \right) }
 \Gamma(s) f^*(1-s) \  x^{-s} ds
=   2 \ (F_c f) \left( x \right) $$$$+   {1\over 2 \pi i} \sqrt {{2\over
\pi}} \int_{\sigma}  \frac { \Gamma (1-s) } { \cos \left(\pi
s/2\right) }  f^*( s) \  x^{s-1} ds.\eqno(2.23)$$
Hence making similar calculations, which were done for the latter integral in (2.16), we  find

$${1\over 2 \pi i}  \int_{\sigma}  \frac { \Gamma (1-s) } { \cos \left(\pi s/2\right) }  \  x^{-s} ds
= {2 x^{-2} \over  \pi}  \int_{0}^\infty   \frac{ y \  K_0(y )}{ \sqrt {y^2 + x^{-2} }} dy =
 { 2 x^{-3/2} \over  \pi}  \   S_{-3/2, -1/2} \left({1\over x}\right),\   x >0.$$
Therefore,
$$F(x) =    2 \ (F_c f) \left( x \right) +  \left({ 2\over \pi}\right)^{3/2} \int_{0}^\infty \sqrt{xt}  \
S_{-3/2,- 1/2} (xt) f(t)dt,$$
which coincides with (2.18).   Similarly,
$$G(x)  = {1\over 2 \pi i}  \left({2\over \pi}\right)^{3/2} \int_{\sigma}  \Gamma(s)\Gamma^2(1-s)\left(\cos \left({\pi s\over 2}\right)
+ \sin \left({\pi s\over 2}\right) \right)^2 \cos \left({\pi s\over
2}\right) f^*(s) \  x^{s-1} ds$$
$$=   {1\over 2\pi i} \sqrt {{2\over \pi}} \int_{\sigma}  \frac {1+  \sin \left(\pi s \right) } { \sin \left(\pi s/2 \right) }
 \Gamma(1-s) f^*(s) \  x^{s-1} ds
=   2 \ (F_s f) \left( x \right)$$$$ +   {1\over 2  \pi i} \sqrt
{{2\over \pi}} \int_{\sigma} \frac { \Gamma (s) } { \cos \left(\pi
s/2\right) }  f^*( 1-s) \  x^{-s} ds$$
and we end up with (2.19),  appealing again to the latter integral in (2.16).  Concerning  inverse operator (2.20),
we write, appealing  to (2.23) and formula (1.8) of the inverse Mellin transform
$$f(x)=  {1\over 2\pi i} \sqrt {{\pi\over 2}} \int_{\sigma}  \frac  { \cos \left(\pi s/2 \right) } {(1+  \sin \left(\pi s \right)) \Gamma(1-s) } F^*(1-s) \  x^{-s} ds =  {1\over 2} (F_c F)(x)$$$$ -  {1\over 2\pi i} {1\over 2\sqrt {2 \pi}} \int_{\sigma}  \frac  { \Gamma(s) \cos \left(\pi s/2 \right) } { \sin^2 \left(\pi (s+1/2)/2  \right) } F^*(1-s) \  x^{-s} ds,\quad  x > 0.\eqno(2.24)$$
In the meantime, the integral
$${1\over 2\pi i}  \int_{\sigma}  \frac  { \Gamma(s) \cos \left(\pi s/2 \right) } { \sin^2 \left(\pi (s+1/2)/2  \right) } \  x^{-s} ds$$
can be calculated by the residue theorem.  It has

$${1\over 2\pi i}  \int_{\sigma}  \frac  { \Gamma(s) \cos \left(\pi s/2 \right) } {\sin^2 \left(\pi (s+1/2)/2  \right)  } \  x^{-s} ds =
 \sum_{k=0}^\infty  {\rm Res }_{ s= - k} \left[  \frac  {\Gamma(s) \cos \left(\pi s/2 \right) } { \sin^2 \left(\pi (s+1/2)/2  \right) } \  x^{-s}\right]$$$$ +   \sum_{k=0}^\infty  {\rm Res }_{ s= - 1/2- 2k} \left[  \frac  {\Gamma(s) \cos \left(\pi s/2 \right) } { \sin^2 \left(\pi (s+1/2)/2  \right) } \  x^{-s}\right],\quad  x >0.\eqno(2.25)$$

The first sum in the right-hand side of the latter equality contains residues in simple poles $s=- k, \ k \in \mathbb{N}_0$ of the gamma-function and by straightforward calculations it gives
$$ \sum_{k=0}^\infty  {\rm Res }_{ s= - k} \left[  \frac  {\Gamma(s) \cos \left(\pi s/2 \right) } { \sin^2 \left(\pi (s+1/2)/2  \right) } \  x^{-s}\right] =  2 \sum_{k=0}^\infty   \frac  {(-1)^k  x^k} { k!}   \cos \left(\pi k/2 \right) = 2 \cos x .$$
The second sum involves double poles $s= -1/2 -2k, \  k \in \mathbb{N}_0$ of the integrand  and we find
$$ {\rm Res }_{ s= - 1/2- 2k} \left[  \frac  {\Gamma(s) \cos \left(\pi s/2 \right) } { \sin^2 \left(\pi (s+1/2)/2  \right) } \  x^{-s}\right]
=  \lim_{ s \to - 1/2- 2k } {d\over ds } \left[ (s+ 1/2+  2k)^2 \frac  { x^{-s}\ \Gamma(s) \cos \left(\pi s/2 \right) } { \sin^2 \left(\pi (s+1/2)/2  \right) } \right]$$
$$= { \sqrt 2(-1)^k \over \pi} x^{2k+1/2} \Gamma(-1/2-2k) \left[  {2\over \pi} \psi( -1/2- 2k)  -  {2\over \pi} \log x   + 1 \right],$$
where $\psi(x)$ is the psi-function (the logarithmic derivative of the gamma-function).  Therefore, substituting these values in (2.25), we obtain
$${1\over 2\pi i}  \int_{\sigma}  \frac  { \Gamma(s) \cos \left(\pi s/2 \right) } {\sin^2 \left(\pi (s+1/2)/2  \right)  } \  x^{-s} ds =
2  \cos x  -   2 \sqrt {{2 x\over  \pi}}  \sum_{k=0}^\infty  {   (-1)^k \  x^{2k} \over (3/2)_ {2k}} $$$$\times  \left[  {2\over \pi} \psi( -1/2- 2k)  -  {2\over \pi} \log x   + 1 \right],\quad  x >0.$$
Hence, returning to (2.24) and employing the generalized Parseval equality (1.9), we come out with inversion formula (2.20). Analogously, we establish (2.21).   The norm inequalities (2.22) are immediate consequences of (2.5)
and isometry property of the Fourier cosine and sine transforms in $L_2$.  To end the proof,  we extend  $F, G $ on the whole  $ L_2(\mathbb{R}_+)$ as  invertible continuous mappings.
\end{proof}

Finally in this section we prove the Plancherel theorem for composition $ \mathcal{H}^2_+ F_c F_s$.   We have

{\bf Theorem 7}.  {\it  The composition $F(x)= \left(\mathcal{H}^2_+ F_c F_s  f\right)(x) $  extends to a bounded invertible map  $F:  L_2(\mathbb{R}_+) \to L_2(\mathbb{R}_+)$ and
$$F (x) = {2^{3/2}  \over \sqrt \pi}  PV  \int_{0}^\infty  \left[ {1\over \pi }  \log \left({x\over t}\right) -1 \right] \frac{ tf(t)  } { x^2-  t^2 }  dt,\quad  x > 0. \eqno(2.26)$$
The inverse operator is written in the form of the integral

$$f(x)=   \sqrt{{2\over \pi^3}} \   PV \int_0^\infty  \  \frac{\sqrt{xt}} {x^2-t^2}   F(t)dt.\eqno(2.27)$$
Moreover,  the composition $F$ satisfies the  norm inequalities $(2.5)$}.

\begin{proof}  Let  $f \in C^{(2)}_c(\mathbb{R}_+)$. Then Theorems 1, 2  yield
$$F(x) = {1\over 2 \pi i}  \left({2\over \pi}\right)^{3/2}  \int_{\sigma}  \left[ \Gamma(s)\Gamma(1-s)\left(\cos \left({\pi s\over 2}\right) + \sin \left({\pi s\over 2}\right) \right)  \cos \left({\pi s\over 2}\right)\right]^2  f^*(s) \  x^{-s} ds$$
$$=  {1\over 2 \pi i} \sqrt {\pi\over 2}   \int_{\sigma} \frac{1+\sin(\pi s)} {\sin^2 \left({\pi s/ 2}\right)}   f^*(s) \  x^{-s} ds=
\left({2 \over \pi}\right)^{3/2}   \int_{0}^\infty  \frac{ \log x-  \log t } { x^2-  t^2 }  tf(t) dt $$
$$- {2\sqrt 2\over \sqrt \pi} PV  \int_{0}^\infty  \frac{ t f(t) } { x^2-  t^2 }  dt.$$
Hence we arrive at (2.26).  Further, to derive (2.27), we have,  reciprocally,

$$f(x)=  {1\over 2 \pi i}  {1\over \sqrt{2\pi} }  {d\over dx }  \int_{\sigma} \frac {\sin^2 \left({\pi s/ 2}\right)\    F^*(s)} { \sin^2 \left(\pi (s+1/2)/2  \right)}   \  {x^{1-s} \over 1-s} ds,\eqno(2.28)$$
where the differentiation is allowed under the integral sign via the absolute and uniform convergence.   Meanwhile, calculating the convergent integral
$${1\over 2 \pi i}  {1\over \sqrt{2\pi} }   \int_{\sigma} \frac {\sin^2 \left({\pi s/ 2}\right)} { \sin^2 \left(\pi (s+1/2)/2  \right)}   \  {x^{1-s} \over 1-s} ds$$
with the use of the residue theorem, involving the left-hand double poles  of the integrand $s= -2k - 1/2, \ k  \in \mathbb{N}_0$,  when $0< x < 1$, we obtain
$${1\over 2 \pi i}  {1\over \sqrt{2\pi} }  \int_{\sigma} \frac {\sin^2 \left({\pi s/ 2}\right)} { \sin^2 \left(\pi (s+1/2)/2  \right)}   \  {x^{1-s} \over 1-s} ds=  - {\sqrt 2\over \pi\sqrt \pi }\left[ 1+ {\log x\over \pi } \right] \sum_{k=0}^\infty {x^{2k+3/2} \over 2k+ 3/2} $$
$$+  {\sqrt 2\over \pi^2 \sqrt \pi } \sum_{k=0}^\infty {x^{2k+3/2} \over (2k+ 3/2)^2 }=
- {\sqrt 2\over \pi\sqrt \pi }\left[ 1+ {\log x\over \pi } \right] \int_0^x  { y^{1/2}\ dy \over 1- y^2}  $$
$$+  {\sqrt 2\over \pi^2 \sqrt \pi } \int_0^x  { y^{1/2} \log(x/y) \ dy \over 1- y^2}, \quad 0 < x < 1.$$
When $x >1$,  we should employ  the right-hand double poles $s= 2k - 1/2, \ k  \in \mathbb{N}$ and the simple pole $s=1$.  This gives the value of the integral
$${1\over 2 \pi i}  {1\over \sqrt{2\pi} }  \int_{\sigma} \frac {\sin^2 \left({\pi s/ 2}\right)} { \sin^2 \left(\pi (s+1/2)/2  \right)}   \  {x^{1-s} \over 1-s} ds=  - \sqrt {{2\over \pi}}  -  {\sqrt 2\over \pi\sqrt \pi }\left[ 1+ {\log x\over \pi } \right] \int_x^\infty   { y^{1/2}\ dy \over  y^2-1 }  $$
$$ +  {\sqrt 2\over \pi^2 \sqrt \pi }  \int_0^{1/x}   { y^{- 1/2}  \log (xy)\over  1- y^2 }   dy , \quad  x >  1.$$
Hence returning to (2.28)  and appealing again to the generalized Parseval equality (1.9), we come out with the inversion formula (2.27) after differentiation under the integral sign, which can be motivated similar  to formulas of the Hilbert transform in $L_2$ (see, \cite{tit},  Th. 90).   The norm inequalities (2.5) follow immediately from  the isometry property of the Fourier transforms in $L_2$.  Extending the composition  on the whole  $ L_2(\mathbb{R}_+)$ as an invertible continuous mapping, we   complete the proof.
\end{proof}

\section{Integral and integro-functional equations}

In this section we will apply the Plancherel theorems for the considered  half-Hartley transform (1.5),  its iteration (2.3) and compositions with the Fourier transforms to investigate the uniqueness and universality of the  closed form solutions of certain singular integral and integro-functional equations.    We begin with an immediate  corollary of Theorem 2.

{\bf Corollary 1}. {\it  Let $g \in   L_2(\mathbb{R}_+) $ be  a given function.  The second kind integral equation with the Stieltjes kernel
$$ f(x) +  {1\over \pi} \int_0^\infty \frac {f(t)} {x+t} dt = g(x),\quad  x > 0,\eqno(3.1)$$
has a  unique solution in $L_2(\mathbb{R}_+) $ given by the formula
$$f(x)=   g (x) -    {2 \over \pi^2}  \int_{0}^\infty  \frac{ \sqrt{xt} \  [ \log x-  \log t] } { x^2-  t^2 }  g(t)  dt.\eqno(3.2)$$
Conversely,   for a given $f \in  L_2(\mathbb{R}_+) $, integral equation $(3.2)$ has a unique solution $g  \in  L_2(\mathbb{R}_+) $ via formula $(3.1)$.}

On the other hand, Theorem 1 leads us the solvability criterium in $L_2(\mathbb{R}_+) $ of the following integro-functional equations  with the Hilbert kernel

$${1\over x} f\left({1\over x}\right) =  {2\over \pi}   \int_{0}^\infty {t f(t) \over t^2- x^2}  dt, \quad  x \in \mathbb{R}_+, \eqno(3.3)$$
$${1\over x} f\left({1\over x}\right)  =  {2\over \pi}   \int_{0}^\infty {x \  f(t) \over x^2- t^2}  dt, \quad  x \in \mathbb{R}_+ .\eqno(3.4)$$
In fact, substituting in (3.3), (3.4) $x$ instead of $1/x$, we arrive at the corresponding second kind homogeneous singular integral equations
$$ f( x) =  {2\over \pi}   \int_{0}^\infty {xt f(t) \over x^2 t^2-  1}  dt, \quad  x  > 0, \eqno(3.5)$$
$$ f( x)  =  {2\over \pi}   \int_{0}^\infty {f(t) \over 1- x^2t^2}  dt, \quad  x >0 .\eqno(3.6)$$

{\bf Corollary 2.} {\it   In order to an  arbitrary function $f \in L_2(\mathbb{R}_+)$ be a solution of either homogeneous  integro-functional equation $(3.3)$ or second kind integral equation $(3.5)$,   it is necessary and sufficient that $f$ has the form of the integral
$$f(x)= {1\over 2\pi i} \int_\sigma \frac{ \varphi (s)}{\cos(\pi s /2)} x^{-s} ds,  \quad   x> 0\eqno(3.7)$$
which is convergent in the mean square sense. It is written  in terms of some function $\varphi(s)$, satisfying  condition $\varphi(s)= \varphi (1-s),\ s \in \sigma$, i.e.  $\varphi (1/2+i\tau)$ is even with respect to $\tau \in \mathbb{R}$.  Analogously,  in order to an  arbitrary function $f \in L_2(\mathbb{R}_+)$ be a solution of either homogeneous  integro-functional equation $(3.4)$ or second kind integral equation $(3.6)$,   it is necessary and sufficient that $f$ has the form of the integral
$$f(x)= {1\over 2\pi i} \int_\sigma \frac{ \rho (s)}{\sin (\pi s /2)} x^{-s} ds,  \quad   x> 0\eqno(3.8)$$
which is convergent in the mean square sense and written  in terms of some function $\rho(s)$, which satisfies   condition $\rho(s)= \rho (1-s),\ s \in \sigma$.}

\begin{proof}  {\it Necessity}.  Let $f \in L_2(\mathbb{R}_+)$  be a solution of  equation (3.3).    In terms of the Mellin transform it can be written as the following functional  equation (see the proof of Theorem 1)
$$ f^*(s) \cot \left({\pi s\over 2}\right) =   f^*(1- s), \quad s \in \sigma.$$
Hence
$$ f^*(s) \cos \left({\pi s\over 2}\right) =   f^*(1- s)\sin \left({\pi s\over 2}\right)= \varphi (s), \quad s \in \sigma\eqno(3.9)$$
and we observe that $\varphi(s)= \varphi (1-s),\ s \in \sigma$.   Therefore,
$$  f^*(s)  = \frac{ \varphi (s)}{\cos(\pi s /2)}$$
and inverting the Mellin transform,  we end up with (3.7).

{\it Sufficiency}.   Conversely,  if $\varphi(1/2+ i\tau)$ is an even function, then from (3.7) we get equalities  (3.9).
Hence the uniqueness theorem for the Mellin transform in $L_2$ drives us at (3.3). The same  concerns to the integral equation (3.5) by virtue of its equivalence to (3.3).  In a similar manner we treat the pair of equations (3.4),  (3.6).
\end{proof}

Theorems 3,4 drive us to the following results.

{\bf Corollary 3.}  {\it  Let $g \in   L_2(\mathbb{R}_+) $ be  a given function.  The second kind integral equations with the Hilbert  kernel
$$ f(x) +  {2\over \pi }   \int_{0}^\infty {x f(t) \over x^2- t^2}  dt = g(x),\quad  x > 0,$$
$$f(x) +  {2\over \pi }   \int_{0}^\infty {t f(t) \over t^2- x^2}  dt= g(x),\quad  x > 0, $$

have    unique solutions  in $L_2(\mathbb{R}_+) $ given by  formulas, respectively,  }
$$f(x)=   {1\over \pi} \int_0^\infty  {\sqrt { x t} \over t^2-  x^2 }  g(t)  dt, $$
$$f(x)=   {1\over \pi}   \int_0^\infty  {\sqrt { x t} \over x^2-  t^2 }  g(t)  dt .$$

{\bf Theorem 8.} {\it  Let $\lambda \in \mathbb{C},\  |1-\lambda|  \neq 1.$ In order to an  arbitrary function $f \in L_2(\mathbb{R}_+)$ be a solution of the  homogeneous  integro-functional equation
 $$f(x) +  {2\over \pi }   \int_{0}^\infty {x f(t) \over x^2- t^2}  dt =   {\lambda \over x} f\left({1\over x} \right),\quad  x > 0, \eqno(3.10)$$
 it is necessary that $f$ has the form of the mean square sense convergent integral
$$f(x)= {1\over 2\pi i} \int_\sigma \frac{ \varphi (s) \   x^{-s}  }{\tan \left({\pi s/ 2}\right) + 1 -\lambda} ds\eqno(3.11)$$
of some function $\varphi(s) \in L_2(\sigma)$,  which satisfies the condition $\varphi(s)= - \varphi (1-s),\ s \in \sigma$.
This condition and the form of solutions $(3.11)$ are also sufficient for those $\varphi$,  whose reciprocal inverse Mellin transform  $\mu(x)$ is a solution of  the integral equation
$$(2- \lambda^2) \mu (x)+  {2\over \pi} \int_0^\infty {\mu (t)\over x+t} dt= 0,  \  x \in  \mathbb{R}_+,\eqno(3.12)$$
where integral $(3.12)$ converges absolutely.   Analogously,  in order to an  arbitrary function $f \in L_2(\mathbb{R}_+)$ be a solution of the  homogeneous  integro-functional equation
$$f(x)   +  {2\over \pi }   \int_{0}^\infty {t f(t) \over t^2- x^2}  dt=  {\lambda\over x} f\left({1\over x} \right),\quad  x > 0, \eqno(3.13)$$
it is necessary  that $f$ has the form of the integral
$$f(x)= {1\over 2\pi i} \int_\sigma \frac{ \varphi (s) \   x^{-s}  }{\cot \left({\pi s/ 2}\right) + 1 -\lambda} ds,\eqno(3.14)$$
which converges in the mean square sense and depends on some function $\varphi(s) \in L_2(\sigma)$,  satisfying  the condition  $\varphi(s)= - \varphi (1-s)$.  This condition and the form of solutions $(3.14)$ are  sufficient for those $\varphi$,  whose reciprocal inverse Mellin transform  $\mu(x)$ is a solution of   integral equation $(3.12)$.}

\begin{proof}    Let $f \in L_2(\mathbb{R}_+)$  be a solution of  equation (3.10).    In terms of the Mellin transform it can be written as the following functional  equation (see the proof of Theorem 3)
$$ f^*(s) \left(1+ \tan \left({\pi s\over 2}\right) \right)=   \lambda f^*(1- s), \quad s \in \sigma.\eqno(3.15)$$
Hence
$$ f^*(s) \tan \left({\pi s\over 2}\right) =  \lambda  f^*(1- s)-  f^*(s)$$
and changing $s$ on $1-s$ in the previous equality, we get
$$ f^*(1-s) \cot  \left({\pi s\over 2} \right)=  \lambda  f^*( s)-  f^*(1-s).$$
Thus, adding these two   equations, we find
$$ f^*(s) \left[ \tan \left({\pi s\over 2}\right) + 1-\lambda \right] +   f^*(1-s) \left[ \cot  \left({\pi s\over 2} \right) + 1 -\lambda\right] = 0.$$
Denoting by $\varphi(s)=  f^*(s) \left[ \tan \left({\pi s/ 2}\right) + 1 -\lambda\right] $, we observe that $\varphi(s)= - \varphi (1-s),\ s \in \sigma$ and $\varphi(s) \in L_2(\sigma)$ if and only if $f^*(s) \in L_2(\sigma)$  because
$$ 0 <\left|  1- | 1-\lambda| \right| \le | \tan \left({\pi s/ 2}\right)+1 -\lambda | \le 2+  |\lambda|, \    s \in \sigma.\eqno(3.16)$$
 Hence $ f^*(s) = \varphi(s) \left[ \tan \left({\pi s/ 2}\right) + 1 -\lambda\right] ^{-1}$ and formula  (1.8)  drives us at  solution (3.11).

 Assuming now the existence of such a function $\varphi (s) \in L_2(\sigma)$ under condition $\varphi(s)= - \varphi (1-s)$, we substitute the value $ f^*(s) = \varphi(s) \left[ \tan \left({\pi s/ 2}\right) + 1 -\lambda\right] ^{-1}$ into equation (3.15).  We have
 $$ \varphi(s) \left[ \frac {1+ \tan \left({\pi s/ 2}\right) }{\tan \left({\pi s/ 2}\right) + 1 -\lambda} + \frac{\lambda}
 {\cot \left({\pi s/ 2}\right) + 1 -\lambda}\right]= 0,$$
 or, via (3.16) and after simple calculations
 $$\varphi(s) \left[ 2 - \lambda^2 + {2\over \sin(\pi s) }  \right]= 0, \quad s \in \sigma.\eqno(3.17)$$
 Taking the inverse Mellin transform of both sides of the latter equality, we arrive at the equation (3.12).  Thus $f(x)$ by formula (3.11) is a solution of integro-functional equation (3.10) for all $\varphi(s)$ under condition $\varphi(s)= - \varphi (1-s)$ such that its inverse Mellin transform is a solution of integral equation (3.12).  The absolute convergence of the corresponding integral follows from the Schwarz inequality.  In the same manner integro-functional equation (3.13) and its solution (3.14) can be treated.

 \end{proof}

 {\bf Corollary 4}. {\it Let $\lambda \in (-\sqrt 2,  \sqrt 2).$ Then the only trivial solution satisfies  integro-functional equations $(3.10), (3.13)$. }

 \begin{proof} When $\lambda =0$, then  condition on $\lambda$ in Theorem 8 fails. However, the solution of (3.10) is trivial via Corollary 3. Otherwise, since $2- \lambda^2 +   2/ \sin(\pi s)  >   0, \ s \in \sigma $, we have from (3.17) $\varphi (s)\equiv 0$ on $\sigma$. Therefore  $f^*(s) \equiv 0$ and the inverse Mellin transform implies  $f=0$, i.e.  the solution of (3.10) is trivial.  The same concerns integro-functional equation  (3.13).
\end{proof}

Composition operator (2.14) is involved to investigate the solvability of the  corresponding homogeneous second kind integral equation
$$\lambda f(x)+  \sqrt { {2\over \pi} } \int_{0}^\infty \left[ \sin(xt) - \cos(xt)+ {2\over \pi}  \sqrt {xt} \   S_{-1/2,1/2} (xt) \right] f(t) dt = 0,\quad  x > 0, \ \lambda \in \mathbb{C}. \eqno(3.18)$$

We have

{\bf Theorem 9.} {\it  Let $|\lambda | <  2.$ In order to an  arbitrary function $f \in L_2(\mathbb{R}_+)$ be a solution of the  homogeneous  integro-functional equation $(3.18)$   it is necessary  that $f$ has the representation
$$f(x)=  {1\over 2\pi i} \int_\sigma \left[   \sqrt {{2\over \pi}}  \Gamma (1-s) \cos \left({\pi s\over 2}\right)\left( 1 +
 \cot\left({\pi s\over 2}\right) \right)- \lambda \right]^{-1}  \varphi (s) x^{-s} ds,  \quad   x> 0,\eqno(3.19)$$
where the integral is convergent in the mean square sense,  depending on some function $\varphi(s) \in L_2(\sigma)$, which satisfies  the condition $\varphi(s)=  \varphi (1-s),\ s \in \sigma$. This condition and the form of solutions $(3.19)$ are  also sufficient for those $\varphi$,  whose reciprocal inverse Mellin transform  $\mu(x)$ is a solution of   integral equation $(3.12)$.}

\begin{proof}   Let $f \in L_2(\mathbb{R}_+)$  be a solution of  (3.18).    Then in terms of the Mellin transform it can be written as follows  (see the proof of Theorem 5)
$$  \sqrt {{2\over \pi}}  \Gamma (s) \sin \left({\pi s\over 2}\right) \left(1 +  \tan \left({\pi s\over 2}\right) \right)   f^*(1-s)  =  - \lambda  f^*( s), \quad s \in \sigma.\eqno(3.20)$$
Hence, changing $s$ on $1-s$ in the previous equation, we find
$$  \sqrt {{2\over \pi}}  \Gamma (1-s) \cos \left({\pi s\over 2}\right) \left(1+  \cot \left({\pi s\over 2}\right) \right)   f^*(s)  =  -\lambda  f^*(1- s).$$
Subtracting one equality from another,   after simple manipulations we end up with
$$  \left[  \sqrt {{2\over \pi}}  \Gamma (1-s)  \cos \left({\pi s\over 2}\right)  \left(1+  \cot \left({\pi s\over 2}\right) \right)
 - \lambda \right]  f^*(s)$$$$
  =  \left[  \sqrt {{2\over \pi}}  \Gamma (s) \sin \left({\pi s\over 2}\right)\left(1+  \tan  \left({\pi s\over 2}\right) \right)
  -\lambda \right] f^*(1- s).$$
Denoting the left-hand side of the previous equation by $\varphi(s)$, we easily verify the condition $\varphi(s) = \varphi(1-s), \ s \in \sigma$.  Moreover, via elementary calculus we derive  $(s= 1/2 +i\tau,\ \tau \in \mathbb{R})$
$$\left|  \sqrt {{2\over \pi}}  \Gamma (1-s) \cos \left({\pi s\over 2}\right) \left(1+   \cot \left({\pi s\over 2}\right) \right)-\lambda  \right| \ge \left|  \sqrt {{2\over \pi}}  \Gamma (1-s) \cos \left({\pi s\over 2}\right) \left(1 +  \cot \left({\pi s\over 2}\right) \right) \right| - |\lambda| $$
$$=   \frac{ 2 \sqrt 2 \cosh (\pi\tau/2)}{\cosh^{1/2} (\pi\tau)} - |\lambda| \ge 2-  |\lambda| > 0$$
and
$$\left|  \sqrt {{2\over \pi}}  \Gamma (1-s) \cos \left({\pi s\over 2}\right)\left(1 +  \cot \left({\pi s\over 2}\right) \right)-\lambda  \right| \le    \frac{ 2 \sqrt 2 \cosh (\pi\tau/2)}{\cosh^{1/2} (\pi\tau)} + |\lambda| \le 2 \sqrt 2 + |\lambda|.$$
Therefore   $\varphi (s) \in L_2(\sigma)$.  Hence calling  inversion formula (1.8) of the Mellin transform, we come out with  solution (3.19) of equation (3.18).    Conversely, assuming  the existence of such a function $\varphi (s) \in L_2(\sigma)$ under condition $\varphi(s)= \varphi (1-s)$, we substitute the value

$$ f^*(s) = \left[ \sqrt {{2\over \pi}}  \Gamma (1-s) \cos \left({\pi s\over 2}\right)\left( 1 +
 \cot\left({\pi s\over 2}\right) \right)- \lambda \right]^{-1} \varphi(s)$$
  into equation (3.20). However, after straightforward calculations   it becomes equation (3.17).  Consequently,   under same conclusions as in  Theorem 8, we complete the proof.
  \end{proof}

  Similarly to Corollary  4 we establish

   {\bf Corollary 5}. {\it Let $\lambda \in (-\sqrt 2,  \sqrt 2).$ Then the only trivial solution satisfies   integral equation $(3.18)$. }

Further, integral operators (2.18), (2.19) are employed  to investigate the $L_2$- solvability of the following homogeneous integral equations of the second kind

$$ 2 \sqrt { {2\over \pi} } \int_{0}^\infty \left[  \cos(xt)+ {1\over \pi}
 \sqrt {xt} \   S_{-3/2, -1/2} (xt) \right] f(t) dt = \lambda f\left( x\right),\quad  x > 0, \ \lambda \in \mathbb{C},\eqno(3.21)$$

$$ 2 \sqrt { {2\over \pi} } \int_{0}^\infty \left[  \sin(xt)+ {1\over \pi}  \sqrt {xt} \
 S_{-1/2, 1/2} (xt) \right] f(t) dt = \lambda f\left( x\right),\quad  x > 0, \ \lambda \in \mathbb{C}.\eqno(3.22)
$$

Precisely,  we arrive at

{\bf Theorem 10.} {\it  Let $|\lambda | <  2.$ In order to an  arbitrary function $f \in L_2(\mathbb{R}_+)$ be a solution of the    integral equation  $(3.21)$   it is necessary  that $f$ has the representation
$$f(x)=  {1\over 2\pi i} \int_\sigma \left[   \sqrt {{2\over \pi}}  \frac{1+\sin \left(\pi s\right)} {
 \cos \left(\pi s/ 2\right) } \Gamma (1-s)  + \lambda \right]^{-1}  \varphi (s) x^{-s} ds,  \quad   x> 0,\eqno(3.23)$$
where the integral is convergent in the mean square sense and depends  on some function $\varphi(s) \in L_2(\sigma)$, which satisfies  the condition $\varphi(s)=  \varphi (1-s),\ s \in \sigma$. This condition and the form of solutions $(3.23)$ are  also sufficient for those $\varphi$,  whose reciprocal inverse Mellin transform  $\mu(x)$ is a solution of  the  integral equation
$$(4- \lambda^2) \mu (x)+  {8\over \pi} \int_0^\infty {\mu (t)\over x+t} dt+   {4\over \pi^2} \int_0^\infty {\log(x/t) \mu (t)\over x-t} dt= 0,  \  x \in  \mathbb{R}_+,\eqno(3.24)$$
where the integrals converge absolutely.  Analogously,  in order to an  arbitrary function $f \in L_2(\mathbb{R}_+)$ be a solution of the    integral equation  $(3.22)$   it is necessary  that $f$ has the representation
$$f(x)=  {1\over 2\pi i} \int_\sigma \left[   \sqrt {{2\over \pi}}  \frac{1+\sin \left(\pi s\right)} {
 \sin \left(\pi s/ 2\right) } \Gamma (1-s)  + \lambda \right]^{-1}  \varphi (s) x^{-s} ds,  \quad   x> 0,\eqno(3.25)$$
where the integral is convergent in the mean square sense and depends  on some function $\varphi(s) \in L_2(\sigma)$, which satisfies  the condition $\varphi(s)=  \varphi (1-s),\ s \in \sigma$. This condition and the form of solutions $(3.24)$ are  also sufficient for those $\varphi$,  whose reciprocal inverse Mellin transform  $\mu(x)$ is a solution of   integral equation $(3.24)$.}

\begin{proof}      Let $f \in L_2(\mathbb{R}_+)$  be a solution of  (3.21).    Using the same technique of the Mellin transform and appealing to the proof of Theorem 6,  we obtain
$$   \sqrt {{2\over \pi}} \frac{1+\sin \left(\pi s\right)} { \sin \left(\pi s/ 2\right) }   \Gamma (s) f^*(1-s)  =  \lambda  f^*( s), \quad s \in \sigma.\eqno(3.26)$$
The  change  $s$ on $1-s$ gives
$$    \sqrt {{2\over \pi}}   \frac{1+\sin \left(\pi s\right)} { \cos \left(\pi s/ 2\right) }\Gamma (1-s) f^*(s)  =  \lambda  f^*(1- s).$$
Subtracting one equality from another,   we define the function $\varphi$ as
$$ \varphi(s)= \left[   \sqrt {{2\over \pi}}  \frac{1+\sin \left(\pi s\right)} {
 \cos \left(\pi s/ 2\right) } \Gamma (1-s)  + \lambda \right]  f^*(s)$$$$
  =  \left[   \sqrt {{2\over \pi}}  \frac{1+\sin \left(\pi s\right)} {
 \sin \left(\pi s/ 2\right) } \Gamma (s)  + \lambda \right] f^*(1- s),$$
which evidently satisfies the equation  $\varphi(s) = \varphi(1-s), \ s \in \sigma$.  Moreover, via elementary calculus we derive  $(s= 1/2 +i\tau,\ \tau \in \mathbb{R})$
$$\left| \sqrt {{2\over \pi}}  \frac{1+\sin \left(\pi s\right)} {
 \cos \left(\pi s/ 2\right) } \Gamma (1-s)  + \lambda    \right| \ge \left|   \sqrt {{2\over \pi}}  \frac{1+\sin \left(\pi s\right)} {
 \cos \left(\pi s/ 2\right) } \Gamma (1-s) \right| - |\lambda| $$
$$=   \frac{ 4 \cosh^2 (\pi\tau/2)}{\cosh (\pi\tau)} - |\lambda| \ge 2-  |\lambda| > 0$$
and
$$\left| \sqrt {{2\over \pi}}  \frac{1+\sin \left(\pi s\right)} { \cos \left(\pi s/ 2\right) } \Gamma (1-s)  + \lambda \right| \le    \frac{ 4 \cosh^2 (\pi\tau/2)}{\cosh (\pi\tau)}  + |\lambda| \le 4 + |\lambda|.$$
Therefore   $\varphi (s) \in L_2(\sigma)$.  Hence calling  inversion formula (1.8) of the Mellin transform, we come out with  solution (3.23) of equation (3.21).    Conversely, assuming  the existence of such a function $\varphi (s) \in L_2(\sigma)$ under condition $\varphi(s)= \varphi (1-s)$, we substitute the value

$$ f^*(s) = \left[   \sqrt {{2\over \pi}}  \frac{1+\sin \left(\pi s\right)} {
 \cos \left(\pi s/ 2\right) } \Gamma (1-s)  + \lambda \right]^{-1}  \varphi(s)$$
 into equation (3.26). But after straightforward calculations   it becomes
 $$ \varphi(s)\left[ \lambda^2-4 - {4\over \sin^2(\pi s)} - {8\over \sin(\pi s)} \right] =0, \ s \in \sigma. $$

 Taking the inverse Mellin transform of both sides of the latter equality we derive integral equation (3.24) (cf. (2.6), (2.7)).
  In the same manner, we examine integral equation (3.22) and its solution (3. 25).
 \end{proof}

   {\bf Corollary 6}. {\it Let $\lambda \in (- 2,  2).$ Then the only trivial solution satisfies   integral equations $(3.21),\ (3.22)$. }

   The final result is the solvability of the integro-functional equation, corresponding the composition operator (2.26)
   $$ {2^{3/2}  \over \sqrt \pi}  \int_{0}^\infty  \left[ {1\over \pi }  \log \left({x\over t}\right) -1 \right] \frac{ tf(t)  } { x^2-  t^2 }  dt= {\lambda\over x} f \left({1\over x}\right) ,\quad  x > 0,\ \lambda \in \mathbb{C}. \eqno(3.27)$$

{\bf Theorem 11.} {\it  Let $|\lambda | <  \sqrt{2\pi}.$ In order to an  arbitrary function $f \in L_2(\mathbb{R}_+)$ be a solution of the    integro-functional  equation  $(3.27)$   it is necessary  that $f$ has the representation
$$f(x)=  {1\over 2\pi i} \int_\sigma \left[   \sqrt {{ \pi\over 2}}  \frac{1+\sin \left(\pi s\right)} {
 \sin^2 \left(\pi s/ 2\right) }   + \lambda \right]^{-1}  \varphi (s) x^{-s} ds,  \quad   x> 0,\eqno(3.28)$$
where the integral is convergent in the mean square sense and depends  on some function $\varphi(s) \in L_2(\sigma)$, which satisfies  the condition $\varphi(s)=  \varphi (1-s),\ s \in \sigma$. This condition and the form of solutions $(3.28)$ are  also sufficient for those $\varphi$,  whose reciprocal inverse Mellin transform  $\mu(x)$ is a solution of  the  integral equation
$$(2\pi - \lambda^2) \mu (x)+  4 \int_0^\infty {\mu (t)\over x+t} dt+   {2\over \pi} \int_0^\infty {\log(x/t) \mu (t)\over x-t} dt= 0,  \  x \in  \mathbb{R}_+,\eqno(3.29)$$
where the integrals converge absolutely. }

\begin{proof}      Let $f \in L_2(\mathbb{R}_+)$  be a solution of  (3.28).    Similarly as above (see  the proof of Theorem 7),  we find
$$   \sqrt {{\pi\over 2}} \frac{1+\sin \left(\pi s\right)} { \sin^2 \left(\pi s/ 2\right) } f^*(s)  =  \lambda  f^*(1- s), \quad s \in \sigma.\eqno(3.30)$$
Changing   $s$ on $1-s$, we have
$$  \sqrt {{\pi\over 2}} \frac{1+\sin \left(\pi s\right)} { \cos^2 \left(\pi s/ 2\right) } f^*(1-s)  =  \lambda  f^*( s).$$
Hence,   we define the function $\varphi$ as
$$ \varphi(s)=   \left[   \sqrt {{\pi\over 2}}  \frac{1+\sin \left(\pi s\right)} {
 \sin^2 \left(\pi s/ 2\right) }   + \lambda \right] f^*(s)
  =  \left[   \sqrt {{\pi\over 2}}  \frac{1+\sin \left(\pi s\right)} {
 \cos^2 \left(\pi s/ 2\right) }   + \lambda \right]f^*(1- s),$$
and clearly   $\varphi(s) = \varphi(1-s), \ s \in \sigma$.  Moreover,   $(s= 1/2 +i\tau,\ \tau \in \mathbb{R})$
$$\left| \sqrt {{\pi\over 2}}   \frac{1+\sin \left(\pi s\right)} { \sin^2 \left(\pi s/ 2\right) }   + \lambda    \right| \ge   \frac{2\sqrt{2\pi} \cosh^2 (\pi\tau/2)}{\cosh (\pi\tau)} - |\lambda| \ge \sqrt{2\pi} -  |\lambda| > 0$$
and
$$\left| \sqrt {{\pi\over 2}}   \frac{1+\sin \left(\pi s\right)} { \sin^2 \left(\pi s/ 2\right) }   + \lambda    \right| \le   \frac{2\sqrt{2\pi} \cosh^2 (\pi\tau/2)}{\cosh (\pi\tau)} + |\lambda| \le 2 \sqrt{2\pi} +  |\lambda| > 0.$$
Therefore   $\varphi (s) \in L_2(\sigma)$.  Hence calling  inversion formula (1.8) of the Mellin transform, we come out with  solution (3.28) of integro-functional equation (3.27).    Conversely, assuming  the existence of such a function $\varphi (s) \in L_2(\sigma)$ under condition $\varphi(s)= \varphi (1-s)$, we substitute the value

$$ f^*(s) =   \left[   \sqrt {{\pi\over 2}}  \frac{1+\sin \left(\pi s\right)} {
 \sin^2 \left(\pi s/ 2\right) }   + \lambda \right]^{-1}  \varphi(s)$$
 into equation (3.30). But after straightforward calculations   it becomes
 $$ \varphi(s)\left[ \lambda^2- 2\pi\left[1  + {1 \over \sin^2(\pi s)} + {2 \over \sin(\pi s)}\right] \right] =0, \ s \in \sigma. $$

 Taking the inverse Mellin transform of both sides of the latter equality we come out with  integral equation (3.29) and complete the proof.
 \end{proof}

   {\bf Corollary 8}. {\it Let $\lambda \in (- \sqrt{2\pi},  \sqrt{2\pi} ).$ Then the only trivial solution satisfies   integro-functional  equation $(3.27)$. }

\noindent {{\bf Acknowledgments}}\\
The present investigation was supported, in part,  by the "Centro de
Matem{\'a}tica" of the University of Porto.\\

\end{document}